\documentclass[12pt, a4paper,twoside]{article}
\usepackage[english]{babel}
\usepackage{amssymb,latexsym}
\usepackage{amssymb}
\usepackage{amsmath}
\usepackage{textcomp}
\usepackage{pst-3d}
\usepackage{latexsym}
\usepackage{pb-diagram}
\newtheorem{thm}{Theorem}[section]
\newtheorem{cor}{Corollary}[section]
\newtheorem{lemme}{Lemma}[section]
\newtheorem{rema}{Remark}[section]

\newtheorem{pro}{Proposition}[section]
\def\tr{{\rm {tr}}}

\newcommand{\cercle}{\mathbb{S}}

\newcommand{\C}{\mathbb{C}}

\newcommand{\R}{\mathbb{R}}

\newcommand{\Spin}{\mathrm{Spin}}
\newcommand{\Spinc}{\mathrm{Spin^c}}
\newcommand{\End}{\mathrm{End}}

\newcommand{\iid}{\mathrm{Id}\,}

\newcommand{\Ree}{\mathrm{Re\,}}

\newcommand{\vol}{\mathrm{vol}}
\newcommand{\divv}{\mathrm{div\,}}
\newcommand{\grad}{\mathrm{grad\,}}
\begin{document}

\title{{\bfseries Lower Bounds for the Eigenvalues of the Dirac Operator on $\Spinc$ Manifolds}}

\author{Roger NAKAD}




\maketitle
\begin{center}
Institut \'Elie Cartan, Universit\'e Henri Poincar\'e, Nancy I, B.P 239\\
54506 Vand\oe uvre-L\`es-Nancy Cedex, France.
\end{center}
\begin{center}
 {\bf nakad@iecn.u-nancy.fr}
\end{center}
\vskip 0.5cm
\begin{center}
 {\bf Abstract}
\end{center}
In this paper, we extend the Hijazi inequality, involving the Energy-Momentum tensor, for the eigenvalues of the Dirac operator on $\Spinc$ manifolds without boundary. The limiting case is then studied and an example is given.
\section{Introduction}
On a compact Riemannian spin manifold $(M^n, g)$ of dimension $n\geqslant2$, Th. Friedrich \cite{ferro}
showed that any eigenvalue $\lambda$ of the Dirac operator satisfies
\begin{eqnarray}
 \lambda^2 \geqslant \lambda_1 ^2 := \frac{n}{4(n-1)}\inf_M S_g,
\label{spinfr}
\end{eqnarray}
where $S_g$ denotes the scalar curvature of $M$. The limiting case
of (\ref{spinfr}) is characterized by the existence of a special spinor called real Killing spinor. 
This is a section $\psi$ of the spinor bundle satisfying for every $X\in \Gamma(TM)$,
$$\nabla_X\psi = -\frac{\lambda_1}{n}X\cdot\psi,$$
where $X\cdot\psi$ denotes the Clifford multiplication and $\nabla$ is the spinorial Levi-Civita connection \cite{16}. On the complement set of zeroes of any spinor field $\phi$, we define $\ell^\phi$ the field of symmetric endomorphisms associated with the field of quadratic forms, denoted by $T^\phi$, called the Energy-Momentum tensor which is given, for any vector field $X$, by 
$$T^\phi (X) = g(\ell^\phi (X), X) = \Ree < X\cdot\nabla_X \phi,\frac{\phi}{\vert\phi\vert^2}>.$$
The associated symmetric bilinear form is then given for every $X, Y \in \Gamma(TM)$ by
$$ g(\ell^\phi (X), Y) = \frac 12 \Ree <X\cdot\nabla_Y\phi + Y\cdot\nabla_X\phi,\frac{\phi}{\vert\phi\vert^2}>.$$
Note that if the spinor field $\phi$ is an eigenspinor, C. B\"{a}r showed that the zero set is contained in a countable union of $(n-2)$-dimensional submanifolds and has locally finite $(n-2)$-dimensional Hausdroff density \cite{barrr}. In 1995, O. Hijazi \cite{12} modified the connection $\nabla$ in the direction of the endomorphism $\ell^\psi$ where $\psi$ is an eigenspinor associated with an eigenvalue $\lambda$ of the Dirac operator and established that
\begin{eqnarray}
 \lambda^2 \geqslant \inf_M (\frac 14 S_g + \vert \ell^\psi\vert^2).
\label{oussamascal}
\end{eqnarray}
The limiting case of (\ref{oussamascal}) is characterized by the existence of a spinor field $\psi$ satisfying for all $X \in \Gamma(TM)$,
\begin{eqnarray}
 \nabla_X\psi = -\ell^\psi(X)\cdot\psi.
\label{morel1}
\end{eqnarray}
The trace of $\ell^\psi$ being equal to $\lambda$, Inequality (\ref{oussamascal}) improves Inequality (\ref{spinfr}) since by the Cauchy-Schwarz inequality, $\vert \ell^\psi\vert^2 \geqslant \frac{(tr(\ell^\psi))^2}{n},$ where $tr$ denotes the trace of $\ell^\psi$. N. Ginoux and G. Habib showed in \cite{habib3} that the Heisenberg manifold is a limiting manifold for (\ref{oussamascal}) but equality in (\ref{spinfr}) cannot occur.\\ \\
Using the conformal covariance of the Dirac operator, O. Hijazi \cite{13} showed that, on a compact Riemannian spin manifold $(M^n, g)$ of dimension $n\geqslant3$, any eigenvalue of the Dirac operator satisfies
\begin{eqnarray}
 \lambda^2 \geqslant \frac{n}{4(n-1)}\mu_1,
\label{oussama1}
\end{eqnarray}
where $\mu_1$ is the first eigenvalue of the Yamabe operator given by 
$$L := 4\frac{n-1}{n-2}\bigtriangleup_g + S_g,$$ $\bigtriangleup_g$ is the Laplacian acting on functions. In dimension 2, C. B\"{a}r \cite{bar} proved that any eigenvalue of the Dirac operator on $M$ satisfies
\begin{eqnarray}
 \lambda^2 \geqslant \frac{2\pi\chi(M)}{Area(M, g)},
\label{oussama}
\end{eqnarray}
where $\chi(M)$ is the Euler-Poincar\'e characteristic of $M$. The limiting case of (\ref{oussama1}) and (\ref{oussama}) is also characterized by the existence of a real Killing spinor. In terms of the Energy-Momentum tensor, O. Hijazi \cite{12} proved that, on such manifolds any eigenvalue of the Dirac operator satisfies the following
\begin{eqnarray}
 \lambda^2 \geqslant \left\{
\begin{array}{l}
\ \ \frac 14 \mu_1 +\inf\limits_M \vert \ell^\psi\vert^2\ \ \ \ \ \ \ \ \text{if}\ \ \ n\geqslant3,\\ \\
\frac{\pi\chi(M)}{Area(M, g)}+\inf\limits_M \vert \ell^\psi\vert^2\ \ \ \ \text{if}\ \ \ n=2.
\end{array}
\right.
\label{hij1}
\end{eqnarray}
Again, the trace of $\ell^\psi$ being equal to $\lambda$, Inequality (\ref{hij1}) improves Inequalities (\ref{oussama1}) and (\ref{oussama}). The limiting case of (\ref{hij1}) is characterized by the existence of a spinor field $\overline\varphi$ satisfying for all $X \in \Gamma(TM)$,
\begin{eqnarray}
 \overline\nabla_X\overline\varphi = -\ell^{\overline\varphi}(X)\ \overline\cdot\ \overline\varphi,
\label{morel11}
\end{eqnarray}
where $\overline\varphi = e^{-\frac{n-1}{2}u}\overline\psi$, the spinor field $\psi$ is an eigenspinor associated with the first eigenvalue of the Dirac operator and $\overline \psi$ is the image of $\psi$ under the isometry between the spinor bundles of $(M^n,g)$ and $(M^n,\overline g = e^{2u}g)$.
Suppose that on a spin manifold $M$, there exists a spinor field $\phi$ such that for all $X\in \Gamma(TM),$
\begin{eqnarray}
 \nabla_X\phi = -E(X)\cdot\phi,
\label{morel}
\end{eqnarray}
where $E$ is a symmetric 2-tensor defined on $TM$. It is easy to see that $E$ must be  equal to $\ell^\phi$. If the dimension of $M$ is equal to 2, Th. Friedrich \cite{friedrich} proved that the existence of a pair $(\phi, E)$ satisfying (\ref{morel}) is equivalent to the existence of a local immersion of $M$ into the euclidean space $\R^3$ with Weingarten tensor equal to $E$. In \cite{morelo}, B. Morel showed that if $M^n$ is a hypersurface of a manifold $N$ carrying a parallel spinor, then the Energy-Momentum tensor (associated with the restriction of the parallel spinor) appears, up to a constant, as the second fundamental form of the hypersurface.
G. Habib \cite{habib1} studied Equation (\ref{morel}) for an endomorphism $E$  not necessarily symmetric. He showed that the symmetric part of $E$ is $\ell^\phi$ and the skew-symmetric part of $E$ is $q^\phi$ defined on the complement set of zeroes of $\phi$ by 
$$ g(q^\phi (X), Y) = \frac 12 \Ree <Y\cdot\nabla_X\phi-X\cdot\nabla_Y\phi,\frac{\phi}{\vert\phi\vert^2}>,$$
for all $X, Y \in \Gamma(TM)$. Then he modifies the connection in the direction of $\ell^\psi + q^\psi$ where $\psi$ is an eigenspinor associated with an eigenvalue $\lambda$ and gets that
\begin{eqnarray}
 \lambda^2 \geqslant \inf_M (\frac 14 S_g +\vert \ell^\psi\vert^2 + \vert q^\psi\vert^2).
\label{habib}
\end{eqnarray}
The Heisenberg group and the solvable group are examples of limiting manifolds \cite{habib1}. For a better understanding of the tensor $q^\phi$, he studied Riemannian flows and proved that if the normal bundle carries a parallel spinor, the tensor $q^\phi$ plays the role of the O'Neill tensor of the flow. Here we prove the corresponding inequalities for $\Spinc$ manifolds:
\begin{thm}\label{main1}
 Let $(M^n, g)$ be a compact Riemannian $\Spinc$ manifold of dimension $n\geqslant2$, and denote by $i\Omega$  the curvature form of the connection $A$ on the
 $\cercle^1$-principal fibre bundle $(\cercle^1 M,\pi,M)$. Then any eigenvalue of the Dirac operator to which is attached an eigenspinor $\psi$ satisfies 
\begin{eqnarray}
 \lambda ^2 \geqslant \inf_M\ \Big(\frac 14 S_g -\frac{c_n}{4}\vert \Omega\vert_g + \vert \ell^{\psi}\vert^2+\vert q^{\psi}\vert^2\Big),
\label{rodgy}
\end{eqnarray}
where $c_n =2[\frac n2]^\frac 12$ and $\vert\Omega\vert_g$ is the norm of $\Omega$ with respect to $g$.
\end{thm}
In this paper, we only consider the deformation of the connection in the direction of the symmetric endomorphism $\ell^\phi$ and hence under the same conditions as Theorem \ref{main1}, one gets
\begin{eqnarray}
 \lambda ^2 \geqslant \inf_M\ \Big(\frac 14 S_g -\frac{c_n}{4}\vert \Omega\vert_g + \vert \ell^{\psi}\vert^2\Big).
\label{fin}
\end{eqnarray}
%
In 1999, A. Moroianu and M. Herzlich \cite{8} proved that on $\Spinc$ manifolds of dimension $n\geqslant 3$, any eigenvalue of the Dirac operator satisfies
\begin{eqnarray}\label{herr}
\lambda^2 \geqslant \lambda_1^2 :=\frac{n}{4(n-1)}\mu_1,
\end{eqnarray}
where $\mu _1$ is the first eigenvalue of the perturbed Yamabe
operator defined by $$L^\Omega = L -c_n\vert\Omega\vert_g.$$
The limiting case of (\ref{herr}) is  characterized by the existence of a real Killing spinor $\psi$ satisfying $\Omega\cdot\psi = i\frac{c_n}{2} \vert\Omega\vert_g\psi$. In terms of the Energy-Momentum tensor we prove:
\begin{thm}\label{main}
Under the same conditions as Theorem \ref{main1}, any eigenvalue $\lambda$ of the  Dirac
operator to which is attached an eigenspinor $\psi$ satisfies
\begin{eqnarray}
\label{ana}
\label{inegalite}
\label{dimension2}
 \lambda^2 \geqslant \left\{
\begin{array}{l}
 \frac 14 \mu _1 +\inf_M\,|\ell^\psi |^2\ \ \ \ \ \ \ \text{if}\ \ \ \ n\geqslant3,\\ \\ 
  \frac{\pi \chi(M) }{Area(M, g)}-\frac 12 \frac{\int_M \vert\Omega\vert_g  v_g}{Area(M,g)}+\inf_M\,|\ell^\psi |^2\ \ \ \text{if}\ \ \ n=2,
\end{array}
\right.
\end{eqnarray}
where $\mu _1$ is the first eigenvalue of the perturbed Yamabe
operator.
\end{thm}
 Using the Cauchy-Schwarz inequality in dimension $n\geqslant3$, we have that Inequality (\ref{inegalite}) implies Inequality (\ref{herr}).
As a corollary of  Theorem \ref{main}, we compare the lower bound to a conformal invariant
(the Yamabe number) and to a topological invariant, in case of 4-dimensional manifolds
whose associated line bundle has self dual curvature (see Corollary (\ref{roro}) and Corollary (\ref{roro1})). Finally, we study the limiting case of (\ref{fin}) and (\ref{ana}), and we give an example.\\ \\
Even though the number $\inf_M |\ell^\psi |^2$ is not a nice geometric invariant, it appears naturally in some situations. For example, on hypersurfaces of certain limiting $\Spinc$ manifolds it is easy to see, with the help of the $\Spinc$ Gauss formula, that it is
precisely the second fundamental form. Also, when deforming the Riemannian metric in the direction of the Energy-Momentum tensor, the eigenvalues of the Dirac operator on a $\Spinc$ manifold are then critical (see \cite{2ana}). The author would like to thank Oussama Hijazi for his support and encouragements.
\section{$\Spinc$ geometry and the Dirac operator}
In this section, we briefly introduce basic notions concerning $\Spinc$ manifolds and the Dirac operator. Details can be found in \cite{5}, \cite{16} and 
\cite{19}. \\ \\
Let $(M^n, g)$ be a compact connected oriented Riemannian manifold of dimension $n\geqslant 2$ without boundary. Furthermore, let $SOM$ be the 
$SO_n$-principal bundle over $M$ of positively oriented orthonormal frames. A $\Spinc$ structure of $M$ is a $\Spin_n^c$-principal bundle $(\Spinc M,\pi,M)$
 and a $\cercle^1$-principal bundle $(\cercle^1 M ,\pi,M)$ together with a double covering given by  $\theta: \Spinc M \longrightarrow SOM\times_{M}\cercle^1 M$ such that
$$\theta (ua) = \theta (u)\xi(a),$$
for every $u \in \Spinc M$ and $a \in \Spin_n^c$, where $\xi$ is the $2$-fold covering of $\Spin_n ^c$ over $SO_n\times \cercle^1$. A Riemannian manifold 
that admits a $\Spinc$ structure is called a Riemannian $\Spinc$ manifold.\\ \\
Let $\Sigma^{c} M := \Spinc M \times_{\rho_n} \Sigma_n$ be the associated spinor bundle where $\Sigma_n = \C^{2^{[\frac n2]}}$ and $\rho_n : \Spin_n^c
\longrightarrow  \End(\Sigma_{n})$ the complex spinor representation. A section of $\Sigma^{c} M$ will be called a spinor and the set of all spinors will be 
denoted by $\Gamma(\Sigma^c M)$. The spinor bundle $\Sigma^c M$ is equipped with a natural Hermitian scalar product, denoted by $< . , .>$ and satisfies
$$< X\cdot\psi, \varphi> = - <\psi, X\cdot\varphi> \ \ \text{for every}\ \ X\in \Gamma(TM)\ \ \text{and}\ \ \psi, \varphi\in \Gamma(\Sigma^c M),$$
where $X\cdot\psi$ denotes the Clifford multiplication of $X$ and $\psi$. With this Hermitian scalar product we define an $L^2$-scalar product
$$(\psi, \phi) = \int_M <\psi, \phi> v_g,$$
for any spinors $\psi$ and $\phi$. Additionally, given a connection 1-form $A$ on $\cercle^1 M$, $A: T(\cercle^1 M)\longrightarrow i\R$ and the connection 1-form 
$\omega^M$ on $SO M$ for the Levi-Civita connection $\nabla^M$, induce a connection on the principal bundle $SO M\times_{M} \cercle^1 M$, and hence 
a covariant derivative $\nabla$ on $\Gamma(\Sigma^{c}M)$ \cite{5}, given by 
\begin{equation}
\nabla_{e_i}\psi = 
\Big[b, e_i(\sigma) + 
\frac{1}{4}\sum_{j=1}^n e_j \cdot \nabla^M_{e_i}e_j \cdot \sigma+ \frac 12 A(s_*(e_i))\sigma\Big],
\label{zshg}
\end{equation}
where $\psi = [b,\sigma]$ is a locally defined spinor field, $(e_1,\ldots,e_n)$ is a local oriented orthonormal tangent frame and $s: U\longrightarrow 
\cercle^1 M$ is a local section of $\cercle^1 M$.\\ \\
The curvature of $A$ is an imaginary valued 2-form denoted by $F_A= dA$, i.e., $F_A = i\Omega$, where $\Omega$ is a real valued 2-form on $\cercle^1 M$. We know
 that $\Omega$ can be viewed as a real valued 2-form on $M$ $\cite{5}$. In this case $i\Omega$ is the curvature form of the associated line bundle $L$. It's
 the complex line bundle associated with the $\cercle^1$-principal bundle via the standard representation of the unit circle.
The spinorial curvature $\mathcal R$ associated with the connection $\nabla$, is given by 
$$\mathcal{R}_{X,Y} = \frac{1}{4}\sum_{i,j=1}^ng\big(R_{X,Y}e_i,e_j\big)\ e_i\cdot
e_j\cdot     +\frac{i}{2}\Omega(X,Y).$$
In the $\Spinc$ case, the Ricci identity translates to 
 \begin{eqnarray}
\sum_j e_j\cdot\mathcal {R}_{e_j,X}\psi = \frac 12 Ric(X)\cdot\psi -\frac i2 (X\lrcorner\Omega)\cdot\psi,
\label{bianchi}
\end{eqnarray}
where $\lrcorner$ denotes the interior product. For every spinor $\psi$, the Dirac operator is locally defined 
by   $$D\psi =\sum_{i=1}^n e_i \cdot \nabla_{e_i} \psi.$$
It is an elliptic, self-adjoint operator with respect to the $L^2$-scalar product and verifies the Schr\"{o}dinger-Lichnerowicz formula 
\[
{D}^2 = {\nabla}^*\nabla + \frac 14 S_g\; \iid_{\Gamma (\Sigma^{c} M)}+ \frac{i}{2}\Omega\cdot,
\]
where $\Omega\cdot$ is the extension of the Clifford multiplication to differential forms given by 
$(e_i ^* \wedge e_j ^*)\cdot\psi = e_i\cdot e_j\cdot\psi$.
\section{Eigenvalue estimates on $\Spinc$ manifolds}
In this section, we prove the lower bound (\ref{rodgy}). This proof is based on the following Lemma given by A. Moroianu and M. Herzlich in \cite{8}:
\begin{lemme}$\cite{8}$. Let $(M^n, g)$ be a $\Spinc$ manifold. For any spinor $\psi \in \Gamma(\Sigma^{c}M)$ and a real 2-form $\Omega$, we have
\begin{eqnarray}
 <i\Omega\cdot\psi,\psi>\ \geqslant -\frac{c_n}{2} \vert \Omega \vert_g\vert \psi\vert^2,
\label{asese}
\end{eqnarray}
where $\vert \Omega \vert_g$ is the norm of  $\Omega$, with respect to $g$ given by $\vert \Omega \vert_g^2=\sum_{i<j} (\Omega_{ij} )^2,$
in any orthonormal local frame. Moreover, if equality holds in (\ref{asese}), then 
\begin{eqnarray}
 \Omega\cdot\psi = i \frac {c_n}{2}\vert\Omega\vert_g \psi.
\label{asese=}
\end{eqnarray}
\label{omega1}
\end{lemme}
\vspace{-0.8cm}
{\bf Proof of Theorem \ref{main1}}: Let $E$ (resp. $Q$) be a symmetric (resp. skew-symmetric) $2$-tensor defined on $TM$. For any spinor field $\phi$, the modified connection $$\widetilde\nabla_X\phi := \nabla_X \phi + E(X)\cdot\phi + Q(X)\cdot\phi,$$ satisfies
$\vert\widetilde\nabla \phi\vert^2 = \vert \nabla\phi\vert^2 -\vert E\vert^2 \vert\phi\vert^2-\vert Q\vert^2 \vert\phi\vert^2.$  After integration  on $M$, the Schr\"{o}dinger-Lichnerowicz formula gives
\begin{eqnarray*}
 \int_M\vert\widetilde\nabla\phi\vert^2 v_g
= \int_M\vert D\phi\vert^2 v_g-\int_M\frac14 S_g\vert\phi\vert^2 v_g &-&\int_M(\vert E\vert^2+ \vert Q\vert^2)\vert\phi\vert^2 v_g \\ &-&\int_M <\frac i2\Omega\cdot\phi, \phi> v_g.
\end{eqnarray*}
Let $\psi$ be an eigenspinor corresponding to the eigenvalue $\lambda$ of $D$. For $E = \ell^\psi$, $Q =q^\psi$ and by Lemma \ref{omega1}, it follows
\begin{eqnarray*}
\lambda^2 \int_M\vert\psi\vert^2 v_g &\geqslant& \frac 14 \int_M S_g\vert\psi\vert^2 v_g +\int_M(\vert \ell^{\psi}\vert^2 +  \vert q^{\psi}\vert^2)
\vert\psi\vert^2 v_g \\ &\ & +\int_M <\frac i2\Omega\cdot\psi, \psi> v_g\\
&\geqslant& \int_M\Big(\frac 14 S_g -\frac{c_n}{4}\vert\Omega\vert_g +\vert \ell^{\psi}\vert^2 +  \vert q^{\psi}\vert^2\Big)\vert\psi\vert^2 v_g.
\end{eqnarray*}
Finally,
$$\lambda ^2 \geqslant \inf_M\Big(\frac 14 S_g -\frac{c_n}{4}\vert \Omega\vert_g + \vert \ell^{\psi}\vert^2 +  \vert q^{\psi}\vert^2\Big).$$
\section{Conformal geometry and eigenvalue estimates}
Before proving  Theorem \ref{main}, we give some basic facts on conformal $\Spinc$ geometry.
The conformal class of $g$ is the set of metrics $\overline g=e^{2u}g$, for a real  
function $u$ on $M$. At a given point $x$ of $M$, we consider a
$g$-orthonormal basis $\{e_1,\ldots,e_n\}\;$ of $T_xM$. The corresponding
$\overline g\,$-orthonormal basis is  denoted by 
$\{\overline{e}_1=e^{-u}e_1,\ldots,\overline{e}_n=e^{-u}e_n\}\;$. This correspondence extends to 
the $\Spinc$ level to give an isometry between the corresponding spinor
 bundles. We put a ``\ $^{\overline{\;\;}}$\ '' above every object which is
 naturally associated with the metric $\overline g$, except for the scalar curvature where $S_g$ (resp. $S_u$ or $S_h$) denotes the scalar curvature associated  with the metric $g$ (resp. $\overline g = e^{2u}g = h^{\frac{4}{n-2}}g$). Then, for any 
spinor fields $\psi $ and $\varphi $, one has
$$<\overline{\psi} ,\overline{\varphi} >=<\psi ,\varphi >\,,$$ 
where $<., .>$ denotes the natural Hermitian scalar products on 
$\Gamma (\Sigma^c M)$, and on $\Gamma (\Sigma^c \,\overline M)$.
The corresponding Dirac operators satisfy   
$$\overline D\,(\,e^{-\frac{(n-1)}{2}u}\;\overline \psi \,)=\;e^{-\frac{(n+1)}{2}u}\; \overline {D\psi}.$$
The norm of any real 2-form $\Omega$ with respect to $g$ and $\overline g$ are related by
$$\vert\Omega\vert_{\overline g} \ = e^{-2u}\vert\Omega\vert_g.$$
O. Hijazi \cite {12} showed that on a spin manifold the Energy-Momentum tensor verifies 
$$|\ell^{\overline{\varphi} }|^2=e^{-2u}\,|\ell^\varphi |^2\,=e^{-2u}\,|\ell^\psi |^2,$$
where $\varphi =e^{-\frac{(n-1)}{2}u}\psi$. We extend the result to a $\Spinc$ manifold and get the same relation.
\begin{lemme}\label{chaftar}
Under the same conditions as Theorem \ref{main1}, any eigenvalue $\lambda$ of the Dirac operator to which is attached an eigenspinor $\psi$ satisfies
$$\lambda^2 \geqslant \frac 14 \sup_u \inf_M (S_u e^{2u} - c_n \vert\Omega\vert_g) + \inf_M \vert \ell^\psi\vert^2.$$
\end{lemme}
Proof: For any spinor field $\phi$ and for any symmetric $2$-tensor $E$ defined on $TM$, the modified connection introduced in \cite{12}: $$\nabla^{E}_X \phi = \nabla_X\phi + E(X)\cdot\phi,$$ verifies $\vert\nabla^E\phi\vert^2 = \vert\nabla\phi\vert^2 - \vert E\vert^2\vert\phi\vert^2$. Using the Schr\"{o}dinger-Lichnerowicz formula on $M$, applied to the spinor field $\overline \phi$ with respect to the metric $\overline g$, yields
\begin{eqnarray}
\int_M\vert\overline\nabla^{E} \overline\phi\vert^2 v_{\overline g}
= \int_M\vert \overline D \ \overline \phi \vert^2v_{\overline g}\ &-&\int_M\frac14 S_u \vert\overline\phi\vert^2v_{\overline g}\ -\int_M\vert E\vert^2\vert\overline\phi\vert^2v_{\overline g}\nonumber \\ &-&\int_M <\frac i2\Omega\ \overline \cdot\ \overline\phi, \overline\phi>v_{\overline g}.
\label{christine}
\end{eqnarray}
For the spinor $\varphi 
=e^{-\frac{(n-1)}{2}u}\;\psi\, $ with $D\psi =\lambda \psi $, 
one gets $\overline D\;\overline \varphi=\;\lambda e^{-u}\; \overline \varphi$, and hence by Lemma \ref{omega1} and for $E= \ell^{\overline \varphi}$
\begin{eqnarray}
\int_M\ \Big[\lambda^2 - (\frac 14 S_u e^{2u} +\vert \ell^{\psi}\vert^2 -\frac{c_n}{4}\vert\Omega\vert_{ g})\Big]e^{-2u} \vert\overline\varphi\vert^2v_{\overline g} \geqslant 0.
\label{uff}
\end{eqnarray}
\begin{lemme}
Let $(M^n, g)$ be a compact Riemannian  $\Spinc$ manifold of dimension $n\geqslant 2$ and $S_g$ (resp. $S_u$ or $S_h$) the scalar curvature associated with the metric $g$ (resp. $\overline g = e^{2u}g = h^{\frac{4}{n-2}}g$). The 2-form $i\Omega$ denotes  the curvature form on the $\cercle^1$-principal bundle associated with the $\Spinc$ structure. We have 
\begin{eqnarray}
\sup_u \inf_M (S_u e^{2u} - c_n\vert\Omega\vert_g) = \left\{
\begin{array}{l}
\mu_1\ \ \ \ \ \ \ \ \ \ \ \ \ \ \ \ \ \ \ \text{if}\ \ \ n\geqslant3,\\ \\
\frac{4\pi \chi(M) - 2\int_M \vert\Omega\vert v_g}{Area(M, g)}\ \ \text{if}\ \ \ \ n=2,
\end{array}
\right.
\end{eqnarray}
where $\mu_1$ is the first eigenvalue of the perturbed Yamabe operator $L^\Omega$.
\label{supinf=}
\end{lemme}
Proof: For $n\geqslant3$, let ${\bf h} > 0$ be an eigenfunction of $L^\Omega$ associated with the eigenvalue $\mu_1$ such that $\int_M {\bf h}^2 v_g =1$. For a conformal metric $\overline g  =e^{2u} g = h^{\frac{4}{n-2}}g$, we have $$S_h h^{\frac{4}{n-2}}-c_n \vert\Omega\vert_g = S_u e^{2u}-c_n \vert\Omega\vert_g =h^{-1}L^{\Omega}h.$$
So $\mu_1 = {\bf h}^{-1}L^{\Omega}{\bf h} = S_{{\bf h}}{\bf h}^{\frac{4}{n-2}}-c_n\vert\Omega\vert_g$. For any positive function $H$, we write $fH={\bf h},$ where $f$ is a positive function, and refering to \cite{9} we get
$$\mu_1 = \int (H^{-1}LH)f^2 H^2\ v_g- c_n \int_M \vert\Omega\vert_g f^2 H^2\ v_g+\int_M H^2 \vert df\vert^2\ v_g.$$
Finally, 
$$\mu_1\geqslant\inf_M(H^{-1} L^{\Omega} H)= \inf_M (S_v e^{2v}-c_n \vert\Omega\vert_g),$$
where $e^{2v} = H^\frac{4}{n-2}$, then $\mu_1 = \sup_u\inf_M (S_u e^{2u} -c_n\vert\Omega\vert_g).$
For $n=2$ and for every $u$ we have $S_ue^{2u} =S_g +2\bigtriangleup_g u.$ The Stokes and Gau$\beta$-Bonnet theorems yield
$$\inf_M(S_ue^{2u} - 2\vert\Omega\vert_g)\leqslant\frac{\int_M \Big(S_u e^{2u}- 2\vert\Omega\vert_g\Big)v_g}{Area(M, g)}= \frac{4\pi\chi(M) - 2\int_M \vert\Omega\vert_g v_g}{Area(M, g)}.$$
Let $u_0$ be a solution of the following equation \cite{0}
\begin{eqnarray}
2\bigtriangleup_g u = \frac{\int_M(S_g- 2\vert\Omega\vert_g) v_g}{Area(M, g)} - S_g +2\vert\Omega\vert_g,
\label{audin}
\end{eqnarray}
hence, $$S_{u_0}e^{2u_0} - 2\vert\Omega\vert_g= 2\bigtriangleup_gu_0 + S_g - 2\vert\Omega\vert_g =  \frac{4\pi\chi(M) - 2\int_M \vert\Omega\vert_g v_g}{Area(M, g)}.$$
{\bf Proof of Theorem \ref{main}}: Combining Lemma \ref{supinf=} and Lemma \ref{chaftar}, Theorem \ref{main} follows.
\begin{rema}
Inequality (\ref{fin}) improves Inequality (\ref{herr}), which itself implies the Friedrich $\Spinc$ inequality given by
\begin{eqnarray}
\lambda^2 \geqslant \frac{n}{4(n-1)}\inf_M(S_g -c_n \vert\Omega\vert_g).
\label{taba}
\label{friedrichspinc}
\end{eqnarray}
Equality holds in (\ref{taba}) if and only if equality holds in (\ref{herr}), i.e., if and only if the eigenspinor $\psi$ associated with the first eigenvalue of $D$ is a real Killing spinor and $\Omega\cdot\psi = i\frac{c_n}{2}\vert\Omega\vert_g \psi.$
\end{rema}
\begin{cor} Any eigenvalue of the Dirac operator on a compact Riemannian $\Spinc$ manifold of dimension $n\geqslant3$, satisfies
$$\lambda ^2 \geqslant \frac{1}{4}\vol(M,g)^{-\frac 2n} \Big(Y(M,[g])-c_n\Vert\Omega\Vert_{\frac n2}\Big)+\mathop {\inf }\limits_M\vert \ell^{\psi}\vert^2,$$
where $Y(M,[g])$ is the Yamabe number given by 
$$Y(M,[g])=\inf_{\eta\neq 0}
 \frac{\int_M 4\frac{n-1}{n-2} \vert d\eta \vert^2 + S_g\eta^2}{\Big(\int_M \vert\eta\vert^{\frac{2n}{n-2}}\Big)^{\frac{n-2}{n}} }.$$
\label{roro}
\end{cor}
\vspace{-0.5cm}
Proof : Using the H\"{o}lder inequality, it follows
$$\mu_1 = \inf_{\eta\neq 0} \frac{\int_M 4\frac{n-1}{n-2} \vert d\eta \vert^2 + (S_g -c_n \vert \Omega\vert_g)\eta^2}{\int_M \eta^2}\geqslant \inf_{\eta\neq 0} \frac{\int_M 4\frac{n-1}{n-2} \vert d\eta \vert^2 + (S_g -c_n \vert \Omega\vert)\eta^2}{\Big(\int_M \vert\eta\vert^{\frac{2n}{n-2}}\Big)^{\frac{n-2}{n}} \vol(M,g)^{\frac 2n}}.$$
Using the H\"{o}lder inequality again, we deduce
$$\mu_1 \ \vol(M,g)^{\frac 2n} \geqslant \inf_{\eta\neq 0}
 \frac{\int_M 4\frac{n-1}{n-2} \vert d\eta \vert^2 + S\eta^2}{\Big(\int_M \vert\eta\vert^{\frac{2n}{n-2}}\Big)^{\frac{n-2}{n}} }-c_n \Big(\int_M \vert \Omega\vert^{\frac n2}\Big)^{\frac 2n}
  = Y(M,[g]) -c_n \Vert\Omega \Vert_{\frac n2}.$$
Finally, replacing in (\ref{inegalite}), we get the result.
\begin{cor}
On a compact 4-dimensional $\Spinc$ manifold with self-dual curvature form $i\Omega$, any eigenvalue of the Dirac operator satisfies 
$$\lambda^2 \geqslant \frac{1}{4}\vol(M,g)^{-\frac 12} \Big(Y(M,[g])-4\pi\sqrt{2}\sqrt{c_1 (L)^2}\Big)+\mathop {\inf }\limits_M\vert \ell^{\psi}\vert^2,$$
where $c_1 (L)$ is the Chern number of the line bundle $L$ associated with the $\Spinc$ structure.
\label{roro1}
\end{cor}
Proof: It follows directly from Corollary \ref{roro} and the fact that if $n=4$ and $\Omega$ self-dual, then 
$\int_M \vert\Omega\vert^2_g v_g= 4\pi^2 c_1 (L)^2$ (see \cite{5}).

\section{Equality case}
In this section, we study the limiting case of (\ref{fin}) and (\ref{ana}). An example is then given.
\begin{pro}\label{prop}
Under the same conditions as Theorem \ref{main1},
\begin{eqnarray}
\text{Equality in (\ref{fin}) holds}\  \Longleftrightarrow 
\left\{
\begin{array}{l}
\nabla_X\psi = -\ell^\psi(X)\cdot\psi,\\
\Omega\cdot\psi = i\frac{c_n}{2}\vert\Omega\vert_g\psi,
\end{array}
\right.
\label{case}
\end{eqnarray}
for any $X \in \Gamma(TM)$ and where $\psi$ is an eigenspinor associated with the first eigenvalue of the Dirac operator.
\end{pro}
Proof: If equality in (\ref{fin}) is achieved, the two conditions follow directly. Now, suppose that $\nabla_X\psi = -\ell^\psi(X)\cdot\psi$ and  $\Omega\cdot\psi = i\frac{c_n}{2}\vert\Omega\vert_g\psi$. The condition $\nabla_X\psi = -\ell^\psi(X)\cdot\psi$ implies that $|\psi |^2$  is constant. Denoting by ${\mathcal R}$ The curvature tensor on the $\Spinc$
bundle associated with the connection $\nabla $, one easily gets the 
following relation 
$$\mathcal R_{X,Y}\,\psi +d \ell^\psi (X,Y)\cdot 
\psi +[\ell^\psi (X),\ell^\psi (Y)]\cdot \psi=0, $$
where $d \ell^\psi$ is a 2-form with values in $ \Gamma(TM) $ given by 
$$ d \ell^\psi(X,Y)= (\nabla_X \ell^\psi) Y - (\nabla_Y \ell^\psi)X.$$
Taking $Y=e_j$ and performing its Clifford multiplication 
by $e_j$ yields by the  Ricci identity (\ref{bianchi}) on a $\Spinc$ manifold
\begin{eqnarray}
-\frac 12 Ric(X)\cdot \psi +\frac i2 (X\lrcorner\Omega)\cdot\psi &+&\sum_j {e_j\cdot d \ell^\psi 
(X,e_j)\cdot\psi}\nonumber \\ &+& \sum_j e_j\cdot [\ell^\psi (X),\ell^\psi (e_j)]\cdot \psi=0.
\label{ricci}
\end{eqnarray}
We then decompose the last two terms in (\ref{ricci}) using that $X\cdot \alpha =X\wedge \alpha -X \lrcorner\alpha$ for any form $\alpha$, it follows
$$\sum_j e_j\cdot d \ell^\psi (X,e_j)\cdot \psi=\sum_j [e_j\wedge d \ell^\psi (X,e_j)]\cdot \psi 
  -[X({\tr}\;\ell^\psi )+\divv \ell^\psi (X)]\psi.$$
$$\sum\limits_j {e_j\cdot [\ell^\psi (X),\ell^\psi (e_j)]\cdot 
\psi}=2\,({\tr }\;\ell^\psi )\,\ell^\psi (X)\cdot \psi -2\sum\limits_j {g(X,
\ell^\psi (e_j))\,\ell^\psi (e_j)\cdot \psi}.$$
Taking the scalar product of (\ref{ricci}) with $\psi $, and after seperating real
and imaginary parts, yields for every vector field $X$ the relation
\begin{eqnarray}
 \Big(X({\tr}\;\ell^\psi )+{\divv}\,\ell^\psi (X)\Big)\vert\psi\vert^2 =\frac i2 <(X\lrcorner\Omega)\cdot\psi,\psi>.\
\label{trace}
\end{eqnarray}
But since Equality (\ref{asese=}) holds we compute
\begin{eqnarray*}
 <(X\lrcorner\Omega)\cdot\psi, \psi> &=& <(X\wedge\Omega)\cdot\psi, \psi> - <X\cdot\Omega\cdot\psi, \psi>\\
&=&  <(X\wedge\Omega)\cdot\psi, \psi>  - i\Big[ \frac n2\Big]^\frac 12 \vert\Omega\vert_g <X\cdot\psi, \psi>.
\end{eqnarray*}
After separating real and imaginary parts, $<(X\lrcorner\Omega)\cdot\psi, \psi>$ must vanish.
Using this and  $\sum_{j=1}^n e_j\cdot(e_j\lrcorner\Omega) = 2\Omega$, Clifford multiplication of  (\ref{ricci}) with 
$e_k$, and for $X=e_k$, gives  
\begin{eqnarray*}
-\frac 12 S_g\psi -i\Omega\cdot\psi=\sum\limits_{k,j} {e_j\cdot (e_k}\wedge d \ell^\psi 
(e_j,e_k))\cdot \psi -2({\tr}\;\ell^\psi )^2\psi +2|\ell^\psi |^2
\psi.
\end{eqnarray*}
An easy computation implies that  $\sum\limits_{k,j} {e_j\cdot (e_k}\wedge d \ell^\psi 
(e_j,e_k))\cdot \psi =0$, hence
\begin{eqnarray}
-\frac 12 S_g + \Big[\frac n2\Big]^\frac 12\vert\Omega\vert_g= -2({\tr}\;\ell^\psi )^2 +2|\ell^\psi |^2,
\label{nouna}
\end{eqnarray}
which implies Equality in (\ref{fin}).
\begin{pro}
 On a compact Riemannian $\Spinc$ manifold 
$(M^n,g)$ of dimension $n\geqslant3$, assume that the first eigenvalue $\lambda_1 $ of the Dirac operator to which 
is attached an eigenspinor $\psi $ satisfies the equality case in (\ref{ana}). Then, $|\ell^\psi |$ is constant and if ${\bf h} > 0$ denotes an eigenfunction of the Yamabe operator corresponding to $\mu_1$, then for any vector field $X$
\begin{eqnarray}
 g(X, \ell^\psi(d{\bf h})-\lambda_1 d{\bf h}) = g(\lambda_1 X -\ell^\psi (X), d{\bf h}) =0.
\label{ba2a}
\end{eqnarray}
\end{pro}
Proof: If $n\geqslant3$ and equality holds in (\ref{ana}), we consider the positive function ${\bf v} > 0$ defined by $e^{2\bf v} = {\bf h}^{\frac{4}{n-2}}$ where ${\bf h}$ is an eigenfunction of the Yamabe operator corresponding to $\mu_1$. Inequality (\ref{uff}) with $u= {\bf v}$ gives $\vert \ell^\psi\vert$ is constant, $ \overline\nabla_X\overline\varphi = -\ell^{\overline\varphi}(X)\ \overline\cdot\ \overline\varphi$ and $\Omega\ \overline \cdot\ \overline\varphi = i \frac{c_n}{2}\vert\Omega\vert_{\overline g}\overline\varphi$. By Proposition \ref{prop}, Equality (\ref{nouna}) and (\ref{trace}) can be considered for the conformal metric $\overline g = e^{2{\bf v}} g= {\bf h}^{\frac{4}{n-2}} g$ to get 
$$({ \tr }\;\ell^{\overline \varphi} )^2:=f^2=\frac 14 S_{{\bf v}} - \frac{c_n}{4}\vert\Omega\vert_{\overline g} +|\ell^{\overline \varphi} |^2,$$
$$\grad f =-\divv \ell^{\overline\varphi}.$$
It is straightforward to see that these two equalities give (\ref{ba2a}).\\ \\
{\bf Example}: If the lower bound (\ref{friedrichspinc}) is achieved, automatically equality holds in (\ref{fin}). Here we will give an example where equality holds in (\ref{fin}) but not in (\ref{friedrichspinc}).\\
Let  $(M^3 ,g)=(S^3,can)$ be endowed with its unique spin structure and consider a real Killing spinor $\psi$ with Killing constant $\frac 12$. As the norm of $\psi$ is constant, we may suppose that $\vert\psi\vert =1$. Let $\xi$ be the Killing vector field on $M$ defined by
$$ig(\xi,X) =<X\cdot\psi,\psi>.$$
In \cite{8}, it is shown that:
\begin{enumerate}
\item $id\xi(X,Y) = -<X\wedge Y \cdot\psi,\psi>\ \ \text{for any }\ \ X, Y \in \Gamma(TM)$.
\item  $d{\vert \xi \vert}^2 = -2d\xi (\xi,.)=-2g(\nabla_{\xi} \xi,.) \simeq -2\nabla_{\xi} \xi =0.$
\item  $\xi\cdot\psi = i\psi\ \text{and} \ \vert\xi\vert =1.$
\item  $\xi\cdot\psi = - e_1 \cdot e_2 \cdot\psi, \ \text{where}\ \{ \xi /  \vert \xi \vert,e_1,e_2\}\ \text{is an oriented local orthonormal frame}$.
\end{enumerate}
Let $h$ be a real constant such that $ h > 1$. We define the metric $g^h$ on $M$, by:
 $$\left\{
\begin{array}{rl}
g^h (\xi,X) &= g(\xi,X)\ \text{pour tout}\  X\in \Gamma(TM),\\
g^h (X,Y) &= h^{-2} g(X,Y)\  \text{pour}\  X,Y \perp \xi.
\end{array}\right.$$
Using the following isomorphism
\begin{eqnarray*}
(TM,g) & \longrightarrow & (TM,g^h) \\
Z & \longrightarrow & Z^h = 
\left\{
\begin{array}{c}
Z\ \  \text{si}\ \  Z=\xi ,\\
hZ\  \ \text{si}\ \  Z \perp \xi,
\end{array}
\right.
\end{eqnarray*}
if $u= \{\xi,e_1,e_2$\} is a positive local g-orthonormal frame defined in a neighborhood $U$ of $x$, then $u^h = \{\xi^h = \xi , e_1 ^h = h e_1 , e_2 ^h = h e_2 \}$ is a positive local $g^h$-orthonormal frame defined in a neighborhood $U$ of $x$.\\There exists an isomorphism of vector bundles (see \cite{8}) given by:
\begin{eqnarray*}
\Sigma_g M & \longrightarrow &\Sigma_{g^h} M \\
\psi=[\tilde u,\phi] & \longrightarrow & \psi^h = [{\tilde u}^h ,\phi],
\end{eqnarray*}
satisfying,
$${<\psi_1 ,\psi_2 >}_{\Sigma_g M } = {<\psi_1 ^h ,\psi_2 ^h>}_{\Sigma_{g^h} M} \ \text{and}\ \ (X\cdot\psi)^h =X^h \cdot\psi^h \ \text{ for any} \ X\in \Gamma(TM).$$
The covariant derivative of the spinor $\psi^h =[\tilde u ^h,\phi]$ is given by (see \cite{8}): 
$$\nabla_{X^h} ^h \psi^h =\frac {h^2}{2}X^h \cdot\psi^h + i((1-h^2)\xi)(X^h)\psi^h.$$
Let $\alpha=(1-h^2)\xi$ be a 1-form on $M$. We may view $i\alpha$ as a connection 1-form on the trivial $\cercle^1$ bundle. Let $L = M \times \C$ be the induced trivial line bundle over $M$. We denote by $\sigma$ the global section of $L$ and by $\nabla^0$ the covariant derivative on $L$ induced by the above connection. It satisfies $$\nabla^0_X \sigma =i\alpha (X) \sigma,\ \text{for any}\ \ X\in \Gamma(TM).$$ On the twisted bundle $\Sigma_{g^h} M \otimes L$, we consider the connection $\overline \nabla =\nabla^h \otimes \nabla^0$ and we calculate 
$$\overline \nabla_{e_1^h} (\psi^h \otimes \sigma)=  \frac{ h^2}{2}e_1^h \cdot(\psi^h \otimes \sigma),$$
$$\overline \nabla_{e_2^h} (\psi^h \otimes \sigma)=  \frac{ h^2}{2}e_2^h\cdot(\psi^h \otimes \sigma),$$
$$\overline \nabla_{\xi} (\psi^h \otimes \sigma)=  (\frac{ -3h^2}{2} + 2)\xi\cdot(\psi^h \otimes \sigma).$$
The spinor $\psi^h\otimes\sigma$ is a section of $\Sigma_{g^h}M\otimes L$, which is, of course, the spinor bundle associated to the $\Spinc$ structure with auxiliary line bundle $L^2$. It is easy to see that $\psi^h\otimes\sigma$ is an eigenspinor associated with the eigenvalue $\frac{h^2}{2} - 2$, and it is clear that $\psi^h\otimes\sigma$ is not a real Killing spinor since $h \neq 1$, so $(M, g^h)$ is not is a limiting manifold for the Friedrich $\Spinc$ inequality. But it is a limiting manifold for the lower bound (\ref{fin}), in fact we will prove that (\ref{case}) holds.\\
The complex 2-form $id\alpha$ is the curvature form associated with the connection $\nabla^0$ on $L$. We have:
$$d\alpha\cdot(\psi^h\otimes\sigma) = (1-h^2)d\xi\cdot (\psi^h\otimes\sigma)=i (h^2-1)h^2 \psi^h\otimes\sigma.$$
The norm of $d\alpha$ with respect to the metric $g^h$ is given by
$$\vert d\alpha\vert^2_{g^h} = (1-h^2)^2 \vert d\xi\vert^2_{g^h} = (1-h^2)^2 (d\xi(e_1^h, e_2^h))^2 = h^4 (1-h^2)^2.$$
Since $h> 1$, $\vert d\alpha \vert_{g^h} = h^2 (h^2-1),$ then the second equation of (\ref{case}) is verified.
Futhermore, it is easy to check that 
$$T^{\psi^h\otimes\sigma} (e_1^h) = T^{\psi^h\otimes\sigma} (e_2^h) = g^h (\ell^{\psi^h\otimes\sigma} (e_1^h), e_1^h) = g^h (\ell^{\psi^h\otimes\sigma} (e_2^h), e_2^h) = -\frac{h^2}{2},$$
$$g^h (\ell^{\psi^h\otimes\sigma} (e_1^h), \xi) = g^h (\ell^{\psi^h\otimes\sigma} (e_2^h), \xi) = g^h (\ell^{\psi^h\otimes\sigma} (e_1^h), e_2^h) =0,$$
$$ T^{\psi^h\otimes\sigma} (\xi) = g^h (\ell^{\psi^h\otimes\sigma} (\xi), \xi) = \frac{3h^2}{2}-2.$$
Finally, it is straightforward to verify that the first equation of (\ref{case}) holds: $$-\ell^{\psi^h\otimes\sigma}(e_1^h)\cdot(\psi^h\otimes\sigma) = \frac{h^2}{2} e_1^h\cdot(\psi^h\otimes\sigma) = \overline \nabla_{e_1^h} (\psi^h \otimes \sigma),$$
$$-\ell^{\psi^h\otimes\sigma}(e_2^h)\cdot(\psi^h\otimes\sigma) = \frac{h^2}{2} e_2^h\cdot(\psi^h\otimes\sigma) = \overline \nabla_{e_2^h} (\psi^h \otimes \sigma),$$
$$-\ell^{\psi^h\otimes\sigma}(\xi)\cdot(\psi^h\otimes\sigma) = (\frac{-3h^2}{2} +2) \xi\cdot(\psi^h\otimes\sigma) = \overline \nabla_{\xi} (\psi^h \otimes \sigma).$$

\end{document}